# PROPAGATION OF OUTLIERS IN MULTIVARIATE DATA

By Fatemah Alqallaf, Stefan Van Aelst,[1] Victor J. Yohai[2]
and Ruben H. Zamar[3]

*Kuwait University, Ghent University, University of Buenos Aires
and University of British Columbia*

We investigate the performance of robust estimates of multivariate location under nonstandard data contamination models such as componentwise outliers (i.e., contamination in each variable is independent from the other variables). This model brings up a possible new source of statistical error that we call "propagation of outliers." This source of error is unusual in the sense that it is generated by the data processing itself and takes place *after* the data has been collected. We define and derive the influence function of robust multivariate location estimates under flexible contamination models and use it to investigate the effect of propagation of outliers. Furthermore, we show that standard high-breakdown affine equivariant estimators propagate outliers and therefore show poor breakdown behavior under componentwise contamination when the dimension $d$ is high.

**1. Introduction.** Most statistical methods are built in the context of a given model and therefore are designed to perform well (e.g., be optimal) for this model. Models are also natural "testing grounds" for statistical procedures and therefore have a profound influence in the way data are processed and analyzed.

Classical models assume data are affected by "normal" noise: small-scale fluctuations arising from measurement errors, item-to-item differences and

Received September 2007; revised December 2007.
[1]Research supported by a grant of the Fund for Scientific Research-Flanders (FWO-Vlaanderen) and by IAP research network Grant nr. P6/03 of the Belgian government (Belgian Science Policy).
[2]Research supported by Grants X-094 from the Universidad de Buenos Aires, PIP 5505 from CONICET, Argentina, PICT 21407 from ANPCyT, Argentina and a PAB grant from ANPCyT.
[3]Research supported by NSERC.

*AMS 2000 subject classifications.* Primary 62F35; secondary 62H12.
*Key words and phrases.* Breakdown point, contamination model, independent contamination, influence function, robustness.

other sources of "well behaved" randomness, for example, Gaussian random variables, Gamma random variables, Poisson processes and other "nice" random disturbances. Contamination models, on the other hand, assume the data may also be affected by abnormal noise: large-scale fluctuations that arise from data contamination, uneven data quality, mixed populations, gross errors, etc. Several contamination models have been proposed in the statistical literature. A nice discussion can be found in Barnet and Lewis (1994).

The best known and most important contamination model is the Tukey–Huber model [Tukey (1962) and Huber (1964)]. This model assumes that, on average, a large fraction $(1 - \epsilon)$ of the data is generated from a classical, normal-error-only model. The remaining data, however, can be affected by abnormal noise. In other words, the Tukey–Huber model assumes a mixture distribution with a fully described dominant component and an unspecified minority component. The mixture fraction $\epsilon$ is a loosely specified nuisance parameter (e.g., $0 \le \epsilon < 0.25$). The goal of a robust statistical analysis is to conduct inference on the dominant part of the mixture, filtering out possible abnormal noise generated by the minority component. The Tukey–Huber contamination model had a profound influence in the general strategy underlying most robust statistical procedures: identify outlying *cases*—those coming from the minority mixture component—and downweight their influence. This model also inspired the definition of key robustness concepts such as influence function, gross-error-sensitivity, maxbias and breakdown point.

The Tukey–Huber contamination model

$$X = (1 - B)Y + BZ,$$

was first introduced in the univariate location-scale setup. The unobservable variables $Y, Z$ and $B$ are independent, $Y \sim F$ [a well-behaved location-scale distribution such as $N(\mu, \sigma^2)$], $Z \sim G$ (an unspecified outlier generating distribution) and $B \sim \text{Binomial}(1, \epsilon)$ (a random contamination indicator). Consequently, the observed variable $X$ has the mixture distribution $(1 - \epsilon)F + \epsilon G$. The model was later extended and used in other settings including regression and multivariate location-scatter models. See, for example, Martin, Yohai and Zamar (1989) and He, Simpson and Portnoy (1990).

The rest of the paper is organized as follows. In Section 2 we introduce a family of contamination models that includes the Tukey–Huber and componentwise contamination models as particular cases. In Section 3 we define and derive the influence function of robust multivariate location estimates under nonstandard contamination models. In Section 4 we discuss propagation of outliers and show that standard high breakdown point (BP) robust estimates propagate outliers. In Section 5 we investigate the breakdown properties under componentwise contamination of robust estimates of multivariate location. Section 6 contains some concluding remarks. Some technical derivations and proofs are given in the Appendix.

**2. Alternative contamination frameworks.** The *multivariate* Tukey–Huber model, where $\mathbf{X}, \mathbf{Y}, \mathbf{Z}$ are $d$-dimensional vectors, may be appropriate for small dimensions but has serious limitations in higher dimensions. A main criticism concerns the assumption that the majority of the cases is free of contamination. Another criticism concerns the downweighting of contaminated cases. When $d$ is large, the fraction of perfectly observed cases can be rather small and the downweighting of an entire case may be inconvenient in the case of "fat and short" data tables where the number of variables (columns) is much larger than the number of cases (rows).

We wish to investigate the robustness properties of classical robust estimates of multivariate location under different contamination models. Suppose that the random vector $\mathbf{Y}$ has density

$$f_{\mathbf{Y}}(\mathbf{y}) = h((\mathbf{y} - \mu_0)' \Sigma_0^{-1} (\mathbf{y} - \mu_0)) \tag{1}$$

and we are interested in estimating the multivariate location vector $\mu_0$. However, we cannot observe $\mathbf{Y}$ directly. Instead, we observe the random vector

$$\mathbf{X} = (\mathbf{I} - \mathbf{B})\mathbf{Y} + \mathbf{B}\mathbf{Z} \tag{2}$$

where $\mathbf{B} = \mathrm{diag}(B_1, B_2, \ldots, B_d)$ is a diagonal matrix, $B_1, B_2, \ldots, B_d$ are Bernoulli random variables with $P(B_i = 1) = \epsilon_i$ and the vector $\mathbf{Z}$ has an arbitrary and unspecified outlier generating distribution.

Note that in principle, the *contamination indicator matrix* $\mathbf{B}$ in model (2) could depend on the vector of uncontaminated observations $\mathbf{Y}$. Likewise, the contamination vector $\mathbf{Z}$ could depend on both the contamination indicator matrix $\mathbf{B}$, and the uncontaminated vector $\mathbf{Y}$. In this paper, however, we restrict attention to the simpler case where $\mathbf{Y}$, $\mathbf{B}$ and $\mathbf{Z}$ are independent.

Different assumptions regarding the joint distribution of $B_1, B_2, \ldots, B_d$ give rise to different contamination models. For example, if $B_1, B_2, \ldots, B_d$ are fully dependent, that is, $P(B_1 = B_2 = \cdots = B_d) = 1$, then model (2) reduces to the classical *fully dependent contamination model* (*FDCM*) which underlies most of the existing robustness theory. An important feature of this model is that the probability of an observation being noncontaminated is $1 - \epsilon$ and so, independently from the dimension, the majority of the cases—rows in the data table—are perfectly observed. Another important feature of this model is that the percentage of contaminated cases is preserved under affine equivariant transformations. Therefore, it is natural that methods designed to perform well under FDCM are affine equivariant and check for the possible existence of a minority of contaminated cases to downweight their influence. Downweighting the influence of suspicious cases is a good strategy when $d$ is relatively small, but becomes less attractive when $d$ is large. For example, downweighting an entire case may be unacceptably wasteful if $d$ is very large and $n$ is relatively small.

Another interesting case is the *fully independent contamination model* (*FICM*) where $B_1, B_2, \ldots, B_d$ are independent. Consider the case $P(B_1 = 1) = \cdots = P(B_d = 1) = \epsilon$, then the probability that a case is perfectly observed under this model is $(1-\epsilon)^d$. Clearly, this probability quickly decreases and goes below the critical value $1/2$ as $d$ increases ($d \geq 14$ when $\epsilon = 0.05$ and $d \geq 69$ when $\epsilon = 0.01$). Another feature of FICM is its lack of affine equivariance. While each column in the data table has on average $(1-\epsilon)100\%$ clean data values, linear combinations of these columns may have a much lower percentage of clean data values. A relevant consequence of this is that in FICM there is a potential to propagate outliers when performing linear operations on the original data. Outlier propagation will be discussed further in Section 4. Intermediate contamination models that fall between FDCM and FICM are briefly discussed in Section 6.

**3. The influence function.** The influence function (IF) is a key robustness tool. It reveals how an estimating functional changes due to an infinitesimal amount of contamination [see Hampel et al. (1986)]. The IF of robust multivariate location estimates has only been defined under the classical FDCM. We wish to extend the definition so that it can be derived under other contamination models.

To fix ideas, we consider the class of $M$-estimates of multivariate location [see, e.g., Tatsuoka and Tyler (2000)] defined as

$$(3) \qquad \mu(H) = \arg\min_{\mathbf{m}} E_H\{\rho[\mathrm{d}^2(\mathbf{X}, \mathbf{m}, \Sigma(H))]\},$$

where

$$\mathrm{d}^2(\mathbf{X}, \mathbf{m}, \Sigma) = (\mathbf{X} - \mathbf{m})'\Sigma^{-1}(\mathbf{X} - \mathbf{m})$$

and $\Sigma(H)$ is a Fisher consistent, preliminary or simultaneous, estimating functional of multivariate scatter. Lemma 3 of Alqallaf et al. (2006) shows that when $\mathbf{X}$ has an elliptical distribution, then $\mu(H)$ is Fisher consistent under mild regularity conditions. Moreover, it is easy to show that $\mu(H)$ satisfies the first order condition:

$$(4) \qquad E_H\{\psi[\mathrm{d}^2(\mathbf{X}, \mu(H), \Sigma(H))](\mathbf{X} - \mu(H))\} = \mathbf{0},$$

where $\psi = \rho'$. Note that equation (4) is satisfied by large classes of estimators for multivariate location such as $M$-estimators [Maronna (1976)], $S$-estimators [Davies (1987), Lopuhaä (1989)], CM-estimators [Kent and Tyler (1996)], MM-estimators [Tatsuoka and Tyler (2000), Tyler (2002)] and $\tau$-estimators [Lopuhaä (1991)].

In order to extend the definition of influence function we must first extend the notion of "point-mass contaminated distribution." Let $G_\epsilon$ be the joint distribution of $(B_1, \ldots, B_d)$, let $\mathbf{z} = (z_1, \ldots, z_d)$ be a given fixed vector in $R^d$ and let $H_0$ be the distribution with density given by (1). Call $H(\epsilon, \mathbf{z})$ the distribution of

$$\mathbf{X} = (\mathbf{I} - \mathbf{B})\mathbf{Y} + \mathbf{B}\mathbf{z},$$

where $\text{diag}(\mathbf{B}) = (B_1, \ldots, B_d) \sim G_\epsilon$ and $\mathbf{Y} \sim H_0$ are independent.

The role played by "point-mass contaminated distributions" in FDCM will be played by $H(\epsilon, \mathbf{z})$ in our more general setup.

The influence function $\text{IF}(\mu, \mathbf{z})$ of the estimating functional $\mu(H)$ given by (3) will be defined and derived for contamination configuration distributions $G_\epsilon$ satisfying (1) $P(B_i = 1) = \epsilon$ $(i = 1, \ldots, d)$; and (2) for any sequence $(j_1, j_2, \ldots, j_d)$ of zeros and ones with $d - k$ zeros and $k$ ones: $P(B_1 = j_1, \ldots, B_d = j_d) = \delta_k(\epsilon)$. These assumptions are clearly satisfied in the case of FDCM and FICM. In FDCM we have $\delta_0(\epsilon) = (1 - \epsilon)$, $\delta_1(\epsilon) = \cdots = \delta_{k-1}(\epsilon) = 0$, and $\delta_k(\epsilon) = \epsilon$. In FICM we have

$$\delta_k(\epsilon) = \binom{d}{k}(1-\epsilon)^{d-k}\epsilon^k, \qquad k = 0, 1, \ldots, d.$$

The (generalized) influence function $\text{IF}(\mu, \mathbf{z})$ is defined as

(5) $$\text{IF}(\mu, \mathbf{z}) = \frac{\partial}{\partial \epsilon}\mu(H(\epsilon, \mathbf{z}))\Big|_{\epsilon=0}.$$

Observe that $H(\epsilon, \mathbf{z})$ and $\text{IF}(\mu, \mathbf{z})$ also depend on $G_\epsilon$ and $H_0$ but, for simplicity, this dependence is not reflected in our notation.

In order to derive $\text{IF}(\mu, \mathbf{z})$, let

(6) $$g(H, \mathbf{m}, \Sigma) = E_H\{\psi(\mathrm{d}^2(\mathbf{X}, \mathbf{m}, \Sigma))(\mathbf{X} - \mathbf{m})\}.$$

From (4) we have that

$$g(H(\epsilon, \mathbf{z}), \mu(H(\epsilon, \mathbf{z})), \Sigma(H(\epsilon, \mathbf{z}))) = \mathbf{0}$$

or, equivalently,

(7) $$\delta_0(\epsilon) g(H_0, \mu(H(\epsilon, \mathbf{z})), \Sigma(H(\epsilon, \mathbf{z}))) \\ + \sum_{k=1}^{d} \delta_k(\epsilon) \sum_{I \in \mathcal{I}_k} g(H(I, \mathbf{z}), \mu(H(\epsilon, \mathbf{z})), \Sigma(H(\epsilon, \mathbf{z}))) = \mathbf{0},$$

where $\mathcal{I}_k = \{I = \{i_1, \ldots, i_k\} : i_1 < \cdots < i_k, 1 \leq k \leq d\}$ and where $H(I, \mathbf{z})$ is the distribution function of $\mathbf{X} = (X_1, \ldots, X_d)$ where $X_i = z_i$ if $i \in I$ and $X_i = Y_i$ if $i \notin I$. In particular, $\mathcal{I}_d = \{1, 2, \ldots, d\}$ and $H(\{1, 2, \ldots, d\}, \mathbf{z}) = \delta_\mathbf{z}$, a point mass distribution at $\mathbf{z}$.

The influence function (5) will now be obtained by differentiating (7) at $\epsilon = 0$. In order to do so we must assume that $\Sigma(H)$ is Fisher consistent at the core model $H_0$ [so $\Sigma(H(0, \mathbf{z})) = \Sigma_0$] and that $\Sigma(H(\epsilon, \mathbf{z}))$ is differentiable with respect to $\epsilon$ at $\epsilon = 0$. When performing the differentiation it is important to notice that

$$(8) \qquad g(H_0, \mu(H(\epsilon, \mathbf{z})), \Sigma(H(\epsilon, \mathbf{z})))|_{\epsilon=0} = g(H_0, \mu_0, \Sigma_0) = \mathbf{0}.$$

Moreover, in the Appendix we show that when $H_0$ is elliptically symmetric,

$$(9) \qquad \frac{\partial}{\partial \epsilon} g(H_0, \mu(H(\epsilon, \mathbf{z})), \Sigma(H(\epsilon, \mathbf{z})))\bigg|_{\epsilon=0} = -A_\psi \, \mathrm{IF}(\mu, \mathbf{z}),$$

where $A_\psi$ is a constant that does not depend on $\mu_0$ and $\Sigma_0$.

Under FDCM, we have $\delta_i(\epsilon) = \delta_i'(\epsilon) = 0$ for $i = 1, \ldots, d-1$, $\delta_d(\epsilon) = \epsilon$ and $\delta_d'(0) = 1$. So, using (7), (8) and (9) we obtain

$$\frac{\partial}{\partial \epsilon} g[H(\epsilon, \mathbf{z}), \mu(H(\epsilon, \mathbf{z})), \Sigma(H(\epsilon, \mathbf{z}))]\bigg|_{\epsilon=0} = -A_\psi \, \mathrm{IF}(\mu, \mathbf{z}) + g(\Delta_\mathbf{z}, \mu_0, \Sigma_0) = \mathbf{0},$$

where $\Delta_\mathbf{z}$ is a point-mass distribution at $\mathbf{z}$. Therefore,

$$(10) \qquad \mathrm{IF}(\mu, \mathbf{z}) = \frac{1}{A_\psi} g(\Delta_\mathbf{z}, \mu_0, \Sigma_0) = \frac{1}{A_\psi} \psi(\mathrm{d}^2(\mathbf{z}, \mu_0, \Sigma_0))(\mathbf{z} - \mu_0).$$

Under FICM we have that $\delta_0(0) = \delta_1'(0) = 1, \delta_1(0) = 0$, and $\delta_i(0) = \delta_i'(0) = 0$ (for $i \geq 2$). So, again using (7), (8) and (9) we have

$$\frac{\partial}{\partial \epsilon} g[H(\epsilon, \mathbf{z}), \mu(H(\epsilon, \mathbf{z})), \Sigma(H(\epsilon, \mathbf{z}))]\bigg|_{\epsilon=0}$$
$$= -A_\psi \, \mathrm{IF}(\mu, \mathbf{z}) + \sum_{k=1}^{d} g(H(I_k, \mathbf{z}), \mu_0, \Sigma_0) = \mathbf{0},$$

where $I_k = \{k\}$. Therefore,

$$(11) \qquad \mathrm{IF}(\mu, \mathbf{z}) = \frac{1}{A_\psi} \sum_{k=1}^{d} g(H(I_k, \mathbf{z}), \mu_0, \Sigma_0).$$

REMARK 1. It is worth noticing that under all the considered models the corresponding influence functions can be interpreted as directional (Gateaux) derivatives. It is well known that in the FDCM case the derivative is in the direction of $\Delta_\mathbf{z}$, a point-mass distribution at $\mathbf{z}$. In the case of FICM the derivative is in the direction of

$$\frac{1}{d} \sum_{k=1}^{d} H(I_k, \mathbf{z}),$$

where $H(I_k, \mathbf{z})$ is the distribution of the random vector $\mathbf{Y} \sim H_0$ with its $k$th component replaced by the constant $z_k$.

*Illustrations.* The effect of an infinitesimal amount of contamination on an estimating functional critically depends on the type of contamination. In the following examples we illustrate some of these differences.

Figure 1 compares the influence functions of an $M$-estimator functional under FDCM and FICM. We consider the case where $\mathbf{Y}$ is bivariate normal with mean zero, variances 1 and correlation $r$. We use the $M$-estimator based on Tukey's bisquare loss function $\rho_c(t) = \min(3t^2/c^2 - 3t^4/c^4 + t^6/c^6, 1)$ with $c^2 = 6$. From Figure 1 we see that the influence functions are fully redescending for FDCM [panel (a)], as is well known. However, for FICM the influence functions are not redescending [panels (b) and (c)]. Therefore, a vanishingly small fraction of large coordinatewise contamination may have a persistent influence on the location $M$-estimate.

Since FICM is not affine equivariant, as discussed further in the next sections, the influence function changes with the amount of correlation $r$. Contrary to the classical contamination case, this change is not just a linear transformation by $\Sigma_0^{1/2}$. This is illustrated in panel (b) ($r=0$) and panel (c) ($r=0.9$) of Figure 1. Note that if $r=0$, then the components are almost solely influenced by contamination in the corresponding component, while in the correlated case they are influenced by contamination in both components. In contrast, the influence function of the coordinatewise $M$-estimator is the same under FDCM and FICM and does not change with correlation. It closely resembles the IF of the bivariate $M$-estimator under FICM with $r=0$ (Figure 1b). As can be expected, the components of the coordinatewise $M$-estimator are only influenced by contamination in the corresponding component, regardless of the value of $r$.

Figure 2 shows the effect of the dimension $d$ on the GES of multivariate location estimators under FDCM and FICM when the core model is multivariate standard normal. We compared the affine equivariant multivariate $S$-location estimator with Tukey loss function and the corresponding coordinatewise $S$-estimator. Under FDCM the truncation parameter in the Tukey bisquare loss function determines the breakdown point of the $S$-estimator [see, e.g., Lopuhaä (1989)]. For the multivariate location $S$-estimator we have chosen the value of the truncation parameter that yields a 50% breakdown point under FDCM. Similarly, for the coordinatewise $S$-estimator we have selected the value of the truncation parameter that yields a 50% breakdown point for each coordinate separately. Figure 2 clearly shows that FDCM severely underestimates the maximal influence of a vanishingly small fraction of contamination on the multivariate $S$-location estimator when $d$ is large. Note that the coordinatewise estimator has the same GES under FDCM and FICM and considerably smaller GES than its multivariate counterpart under FICM.

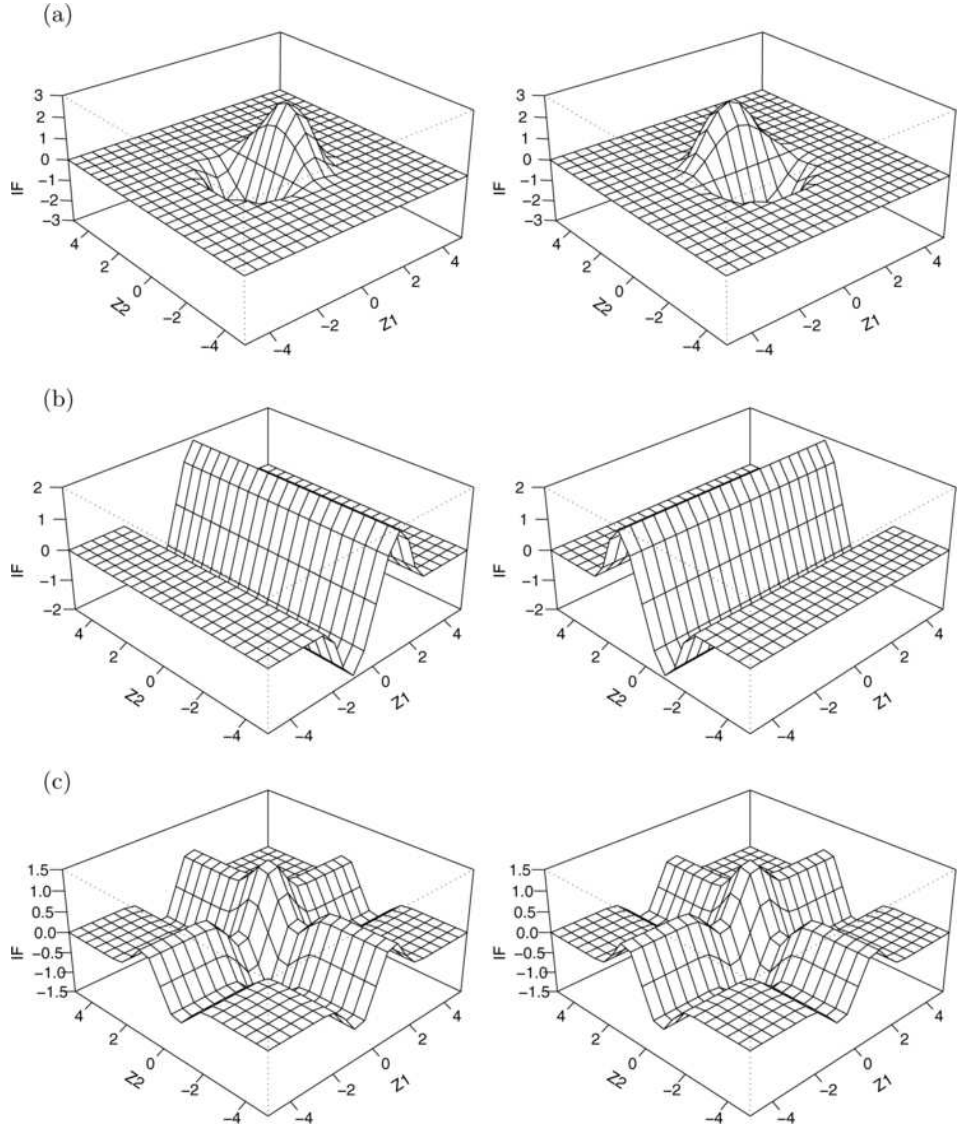

FIG. 1. *Influence functions for Tukey bisquare M-estimator of bivariate location. Left panels are for the first component, right panels are for the second component. Panel* (a) *FDCM; Panel* (b) *FICM with $r = 0$; Panel* (c) *FICM with $r = 0.9$.*

**4. Propagation of outliers.** FDCM is translation-scale equivariant and affine equivariant. Therefore, if a random vector $\mathbf{X}$ follows this model, then an affine transformation $\tilde{\mathbf{X}} = \mathbf{A}\mathbf{X} + \mathbf{b}$ will also follow the model, for any invertible matrix $\mathbf{A}$ and vector $\mathbf{b}$. In particular, if $\mathbf{X}$ has a probability $\epsilon$ of contamination, the same probability holds for $\tilde{\mathbf{X}}$. On the other hand, the independent contamination model is *not affine equivariant*. In fact, suppose that the random vector $\mathbf{X}$ follows the FICM and $\mathbf{A}$ is an invertible $d \times d$ matrix, then the transformed vector

$$\tilde{\mathbf{X}} = \mathbf{A}\mathbf{X} + \mathbf{b} = \mathbf{A}(\mathbf{I} - \mathbf{B})\mathbf{Y} + \mathbf{A}\mathbf{B}\mathbf{Z} + \mathbf{b}$$

is in general different from $(\mathbf{I} - \mathbf{B})\mathbf{A}\mathbf{Y} + \mathbf{B}\mathbf{A}\mathbf{Z} + \mathbf{b}$, unless $\mathbf{A}\mathbf{B} = \mathbf{B}\mathbf{A}$ (i.e., $\mathbf{A}$ is diagonal). Therefore, $\tilde{\mathbf{X}}$ does not follow the independent contamination model.

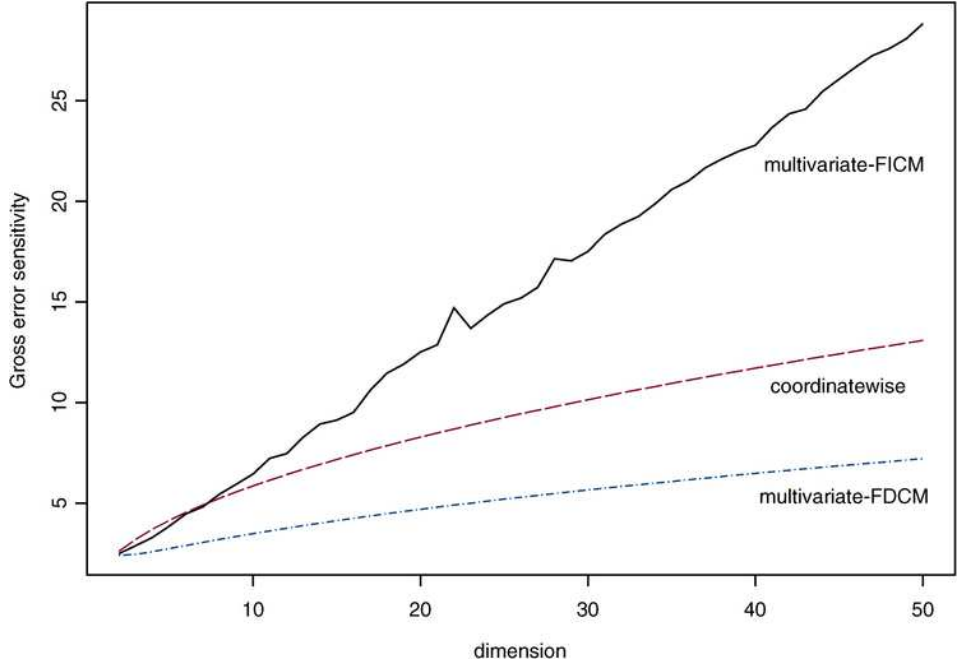

FIG. 2. *Gross error sensitivity for Tukey bisquare* 50% *breakdown S-estimator of multivariate location and corresponding* 50% *breakdown coordinatewise S-estimator.*

The lack of affine equivariance of FICM causes a phenomenon that we call "*outlier propagation*." FICM assumes that each column in the data table contains an average fraction $\epsilon$ of contamination. Since affine transformations linearly combine the columns, the independent contamination property is lost.

To illustrate this, we generated a small two-dimensional data set of size $n = 20$. Both components come from a standard Gaussian distribution and we added independent contamination to each component with a contamination probability of 30%. The contaminated data come from a Gaussian distribution with mean 10 and variance 1. Histograms of the original components $X_1$ and $X_2$ are shown in the top panels of Figure 3. Both histograms show a clear majority of clean data with approximately 1/3 of outlying points on the right. The thick vertical lines indicate the medians, 0.22 for $X_1$ and 0.95 for $X_2$. The medians are slightly affected by the heavy contamination but still summarize well the majority of the data. Now, we consider an affine transformation: $L_1 = 0.64X_1 + 0.77X_2$ and $L_2 = 0.78X_1 + 0.62X_2$. Histograms of the components $L_1$ and $L_2$ are shown in the bottom panels of Figure 3. From these histograms it is clear that both components now contain a majority of contaminated cells and hence do not satisfy FICM with 30% contamination anymore. In fact, we have three distinct groups in each dimension consisting of 49%, 42% and 9% of the data. Note that the medians of $L_1$ and $L_2$ no longer reflect the location of the clean data.

Data following FICM and other nonaffine equivariant versions of model (2) can severely upset standard, affine equivariant robust procedures. To illustrate this, we consider the following example. We generated 100 observations from a 15 dimensional standard Gaussian distribution and added independent contamination to each column with a contamination probability of $\epsilon = 15\%$. The contamination is obtained by adding a constant $t$ to the generated values. The overall probability of a contaminated cell is thus 15%, which is reasonably low, so one might expect to obtain reliable estimates if a robust estimator is used. However, the probability that an observation is contaminated equals $1 - (1-\epsilon)^d > 90\%$. By applying a linear transformation to these data, the outlier propagation effect can spread contamination in one of the components of an observation over all its components. This results in transformed data with a contamination probability of more than 90% for each cell. In such a setting, no robust estimator is supposed to be reliable anymore. However, for affine equivariant robust estimators the original and linearly transformed data sets are equivalent, which has a devastating effect on their performance in this setting, as illustrated in Figure 4. In this figure we generated data sets with 15% of independent contamination in each column, as explained above. We varied the size of the contamination constant $t$ from 0 to 100. We calculated the multivariate location of the data using the sample mean, the coordinatewise median and three affine equivariant robust location estimators: the (i) Minimum Volume Ellipsoid (MVE), (ii) Minimum Covariance Determinant (MCD) [both proposed by Rousseeuw (1984)] and (iii) the Stahel–Donoho estimator, independently proposed by Stahel (1981) and Donoho (1982). If an estimator is not affected by contamination, then all components of the location vector should be close to zero. On the other hand, if the contamination affects the estimator, then some components of the location vector will become biased. For each estimator, we plotted the largest componentwise bias of the estimated location vector against the size of the contamination. Note that the bias of both MVE and MCD increases without bound. The bias of the Stahel–Donoho estimator as implemented in Splus increases even faster and the estimator crashes when the contamination constant exceeds 7. The three affine equivariant robust estimates show clear signs of breaking down. Not surprisingly, so does the sample mean. On the other hand, the coordinatewise median is hardly affected by the outliers in each component. This example clearly shows that robust affine equivariant methods are not robust against propagation of outliers. (A more rigorous treatment of this claim will be given in the next section.) Hence, these methods are not well suited for situations where the contamination regime operates on individual variables (columns) rather than individual cases (rows).

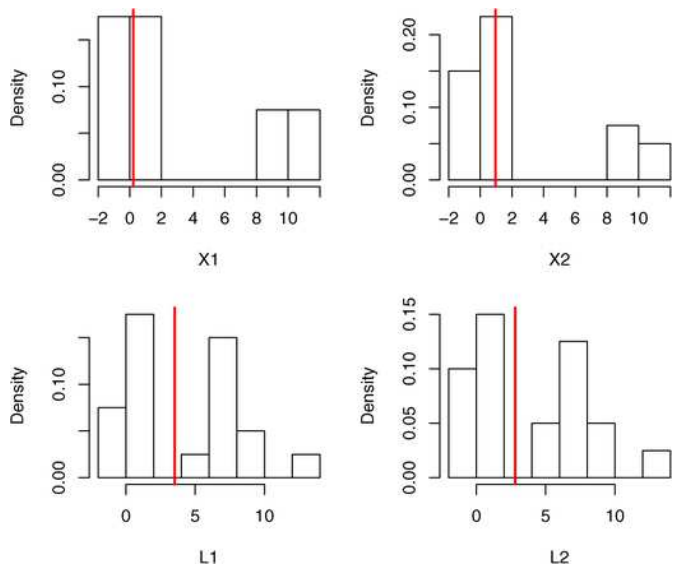

Fig. 3. *FICM outliers propagated by linear combinations.*

**5. Affine equivariance and independent contamination.** For simplicity, we will keep the Section 3 assumption that the marginal probabilities of a contaminated cell are equal for all components, that is, $P(B_1 = 1) = \cdots = P(B_d = 1) = \epsilon$. However, with obvious modifications the results hold for the general case as well.

For each distribution $G_0$ on $R$ with finite first moment, let $\mathcal{G}_h(G_0)$ be the set of distribution functions $G$ on $R^d$ with marginal distributions, which are all stochastically larger than $G_0(x - h)$. For each $\delta > 0$ set

$$\mathcal{F}_{\delta,h}(G_0) = \{H = (1/2 - \delta)H_0 + (1/2 + \delta)G,\ G \in \mathcal{G}_h(G_0)\}.$$

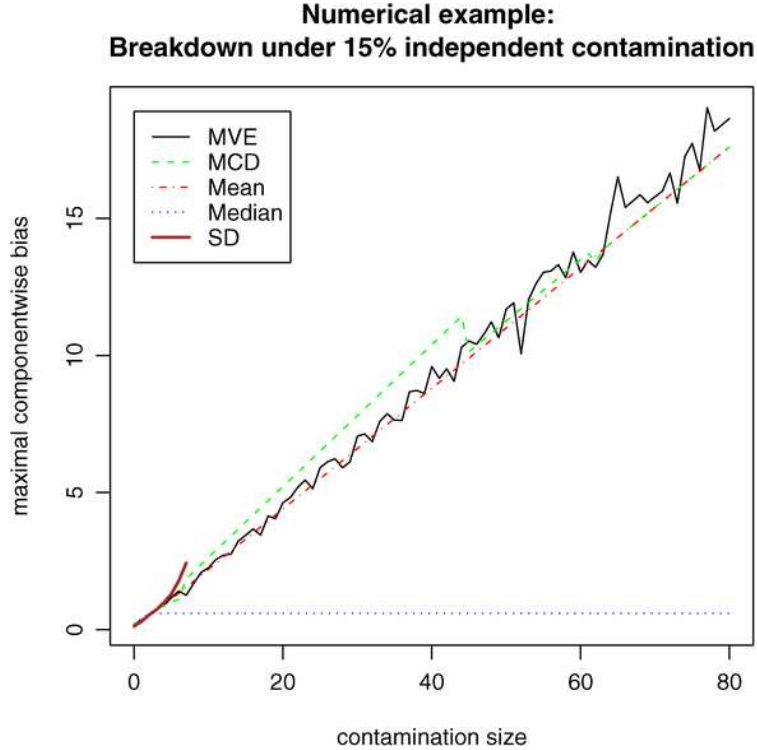

Fig. 4. *Affine equivariant, high breakdown-point estimators try to identify outlying cases and break down when more than 50% of the cases are contaminated, which can easily occur with small fractions of independent contamination in the variables when the dimension is moderately large.*

DEFINITION 1. Let $\mathbf{T} = (T_1, \ldots, T_d)$ be an equivariant multivariate location estimating functional on $R^d$. We say that $\mathbf{T}$ is $\delta$-consistent at infinity, when the central model is $H_0$, if for any distribution $G_0$

$$\lim_{h \to \infty} \inf_{H \in \mathcal{F}_{\delta,h}(G_0)} T_i(H) = +\infty, \qquad 1 \leq i \leq d.$$

In other words, $\delta$-consistent estimates have the property that if at least $1/2 + \delta$ of the mass goes to infinity for all the coordinates, then all the coordinates of the estimate go to infinity too. Note that $\delta_2 > \delta_1$ implies $\mathcal{F}_{\delta_2,h}(G_0) \subset \mathcal{F}_{\delta_1,h}(G_0)$, thus if $\mathbf{T}$ is $\delta_1$-consistent, then it is also $\delta_2$-consistent.

Let us introduce the following notation. Given a distribution $H_0$ on $R^d$ denote by $\mathcal{F}_\epsilon^I$ its FICM contamination neighborhood of size $\epsilon$ that contains all the distributions of $\mathbf{X} = (\mathbf{I} - \mathbf{B})\mathbf{Y} + \mathbf{B}\mathbf{Z}$ where $\mathbf{Y}, \mathbf{B}$ and $\mathbf{Z}$ are independent, $\mathbf{Y}$ has distribution $H_0$, $\mathbf{B}$ is a diagonal matrix where the diagonal elements $B_1, \ldots, B_d$ are independent Bernoulli variables such that $P(B_i = 1) = \epsilon$ and **Z** has an arbitrary distribution $H^*$. We denote by $\mathcal{F}_\epsilon^D$ its FDCM contamination neighborhood that contains the distributions of the form

$$H = (1-\epsilon)H_0 + \epsilon H^*,$$

where $H^*$ is arbitrary.

We can now define the breakdown point under FICM, $\varepsilon^*_{\text{FICM}}$, of a multivariate location estimator $\mathbf{T}(H)$ as the smallest probability $\epsilon$ of contamination in each of the components that is needed to make $\|\mathbf{T}(H)\|$ arbitrary large. That is,

$$\varepsilon^*_{\text{FICM}}(\mathbf{T}, H_0) = \inf\left\{\epsilon > 0;\ \sup_{H \in \mathcal{F}_\varepsilon^I} \|\mathbf{T}(H)\| = +\infty\right\}.$$

Theorem 1 shows that the FICM breakdown point of any equivariant estimate of location which is $\delta$-consistent at infinity under FICM is at most $1 - (1/2 - \delta)^{1/d}$. Hence, if $\delta$ is independent of $d$, the FICM breakdown point tends to 0.

THEOREM 1. *Let $\mathbf{T}(H)$ be an affine equivariant multivariate location estimator that is $\delta$-consistent at infinity for the central distribution $H_0$, with finite first moments. If*

$$\epsilon > \epsilon_0 = 1 - (1/2 - \delta)^{1/d},$$

*then*

$$\sup_{H \in \mathcal{F}_\varepsilon^I} \|\mathbf{T}(H)\| = +\infty.$$

*Hence,*

$$\varepsilon^*_{\text{FICM}}(\mathbf{T}, H_0) \leq 1 - (1/2 - \delta)^{1/d}.$$

PROOF. Consider the linear transformation $\mathbf{U} = A\mathbf{X}$ with

$$A = \begin{pmatrix} 2 & 1 & \cdots & 1 \\ 1 & 2 & \cdots & 1 \\ \vdots & \vdots & & \vdots \\ 1 & 1 & \cdots & 2 \end{pmatrix}.$$

Note that $A$ is invertible since its eigenvalues are 2 with multiplicity one and 1 with multiplicity $(d-1)$.

Let $H_h \in \mathcal{F}_\varepsilon^I$ where $\epsilon > \epsilon_0$ and $\mathbf{Z} \sim \delta_h$ with $\delta_h$ the point mass at $(h, \ldots, h) \in R^d$. It follows that with probability $(1-\epsilon)^d = 1/2 - \delta^*$, with $\delta^* > \delta$ the vector $\mathbf{X}$ comes from $H_0$, and thus with probability $1/2 + \delta^*$ at least one component of $\mathbf{X}$ is equal to $h$.

Let $\tilde{H}_h$ and $\tilde{H}_0$ be the distributions of $\mathbf{U}$ when $\mathbf{X}$ has distribution $H_h$ and $H_0$, respectively. Then

$$\tilde{H}_h = (1 - \delta^*)\tilde{H}_0 + \delta^* G_h,$$

where $G_h$ is the distribution of $\mathbf{U}$ when $\mathbf{X}$ has distribution $H_h$ conditionally on $\sum_{i=1}^d B_i > 0$. Therefore, all the marginals of $G_h$ are stochastically larger than $G_0(u - h)$ where $G_0$ is the distribution of $-2\sum_{j=1}^d |Y_j|$ with $\mathbf{Y} \sim H_0$. Since $\mathbf{T}$ is $\delta^*$-consistent at infinity, we then have

$$\lim_{h \to \infty} \|\mathbf{T}(\tilde{H}_h)\| = +\infty.$$

Since $A$ is invertible and $\mathbf{T}$ is affine equivariant,

$$\lim_{h \to \infty} \|\mathbf{T}(H_h)\| = \lim_{h \to \infty} \|A^{-1}\mathbf{T}(\tilde{H}_h)\| = +\infty,$$

proving the theorem. □

It is obvious that a scatter estimate breaks down whenever the multivariate location estimate it is using to center the data breaks down. Therefore, although Theorem 1 is stated for multivariate location, it has clear implications for the companion scatter estimates.

The following lemma is proven in the Appendix and will be used to show that many of the well-known affine equivariant robust estimators of multivariate location are $\delta$-consistent at infinity.

LEMMA 1. *Suppose that $\mathbf{T}(H)$ is location-scale equivariant and can be represented as a weighted average, that is, it can be written as*

(12) $$\mathbf{T}(H) = E_H(\mathbf{X} w(H, \mathbf{X})),$$

*where the weight function $w(H, \mathbf{x})$ satisfies:* (i) $w(H, \mathbf{x}) \geq 0$, (ii) *there exists $K$ such that $w(H, \mathbf{x}) \leq K$ and* (iii) *there exists $\eta > 0$ such that $P_H(w(H, \mathbf{x}) > \eta) > 1/2 - \delta_0$ for some $\delta_0 > 0$. Then $T$ is $\delta$-consistent at infinity when the central model distribution $H_0$ has finite first moments, for all $\delta > \delta_0$.*

*Examples of $\delta$-consistency at infinity.* The following examples illustrate how Lemma 1 can be used to show $\delta$-consistency at infinity for well-known affine equivariant high-breakdown (under FDCM) estimators of multivariate location. Table 1 shows that for higher dimensions ($d \geq 10$) a small amount of contamination in each variable suffices to break down such estimators.

TABLE 1
*Minimal fraction of independent contamination that causes breakdown of $\delta$-consistent, affine equivariant estimators*

| | **Dimension** | | | | | | | |
|---|---|---|---|---|---|---|---|---|
| | 1 | 2 | 3 | 4 | 5 | 10 | 15 | 20 | 100 |
| $\epsilon$ | 0.50 | 0.29 | 0.21 | 0.16 | 0.13 | 0.07 | 0.05 | 0.03 | 0.01 |

*Coordinatewise mean and median.* It is clear that the sample mean satisfies (12) with weights $w(H, \mathbf{x}) = 1$, hence the assumptions of Lemma 1 hold in this case. Although the coordinatewise median does not satisfy the assumptions of Lemma 1, a simple argument shows that it is $\delta$-consistent at infinity for all $\delta > 0$. Using the notation introduced in the proof of Lemma 1, we have that $P_{G_h}(X_i \leq \sqrt{h}) \to 0$ and so $\lim_{h\to\infty}[(1/2 - \delta)H_{0i}(\sqrt{h}) + (1/2 + \delta)G_{hi}(\sqrt{h})] < 1/2$. Therefore, $\lim_{h\to\infty} \text{Med}(F_{hi}) = \infty$. Note, however, that the coordinatewise median is not affine equivariant and thus Theorem 1 does not apply in this case.

*Minimum Covariance Determinant.* The Minimum Covariance Determinant estimator (MCD) of multivariate location, introduced by Rousseeuw (1984), is defined as a scaled weighted mean $\mathbf{T}_{MCD}(H)$ with weight $w(\mathbf{x}, H) = I_{A^*}(\mathbf{x})/P_H(A^*)$. The set $A^*$ is determined as follows. Let $\mu(H, A) = \int_A \mathbf{x} \, dH(\mathbf{x})/P_H(A)$ be the mean associated to any subset $A \subset R^d$. Then $A^*$ is such that its covariance matrix $\Sigma(H, A^*) = \int_{A^*} (\mathbf{x} - \mu(H, A^*))(\mathbf{x} - \mu(H, A^*))' \, dH(\mathbf{x})$ has smallest determinant among all subsets $A$ such that $P_H(A) \geq 1/2$. Clearly, the weights are nonnegative and bounded. Moreover, since $1 \leq w(\mathbf{x}, H) \leq 2$ for all $\mathbf{x} \in A^*$, we have that $P(w(H, \mathbf{x}) > \eta) > 1/2 - \delta_0$ for any $\eta < 1$ and $\delta_0 > 0$. Then $\mathbf{T}_{MCD}$ satisfies the assumptions of Lemma 1 and is $\delta$-consistent at infinity for any $\delta > 0$.

*S-estimators.* Consider a function $\rho: R \to R^+$ that satisfies the following assumptions:

A1. $\rho$ is even, bounded and nondecreasing on $[0, \infty)$ with $\rho(0) = 0$. Without loss of generality we will take $\rho(\infty) = 1$.
A2. $\rho$ is differentiable, $\psi(t) = \rho'(t)$ is differentiable at 0, and $u(t) = \psi(t)/t$ is nonincreasing on $[0, \infty)$. We will also assume that $\rho(u) < 1$ implies $\psi(u) > 0$.

Then $(\mathbf{T}(H), S(H))$ is defined by the values $(\mu, \Sigma)$ satisfying

(13) $$(\mathbf{T}(H), S(H)) = \underset{\mu, \Sigma}{\arg\min} \det(\Sigma)$$

subject to

(14) $$E_H(\rho(\mathrm{d}(\mathbf{x}, \mu, \Sigma)/s_0) = b.$$

It can be shown [see, e.g., Davies (1987)] that if $\rho$ is differentiable then $\mathbf{T}(H)$ satisfies the following equation

(15) $$\mathbf{T}(H) = E_H(w(\mathbf{X}, H)\mathbf{X})$$

with

(16) $$w(\mathbf{x}, H) = \frac{u(\mathrm{d}(\mathbf{x}, \mathbf{T}(H), S(H)))}{E_H(u(\mathrm{d}(\mathbf{X}, \mathbf{T}(H), S(H))))}$$

and

(17) $$u(\mathbf{x}, H) = \frac{\psi(\mathrm{d}(\mathbf{x}, \mathbf{T}(H), S(H)))}{\mathrm{d}(\mathbf{x}, \mathbf{T}(H), S(H))}.$$

The following lemma (proven in the Appendix) shows that $S$-estimators are $\delta$-consistent at infinity for any $\delta > 0$.

LEMMA 2. *Suppose* A1 *and* A2 *are satisfied. Then the weight function* $w(\mathbf{x}, H)$ *associated with the $S$-location estimate* $\mathbf{T}(H)$ *with* $b = 1/2$ *satisfies assumptions* (i), (ii) *and* (iii) *of Lemma 1, for any* $\delta > 0$.

Lemma 2 can be extended to $\tau$-estimates of multivariate location as defined by Lopuhaä (1991). In addition, a simple argument shows that Rousseeuw's minimum volume ellipsoid is also $\delta$-consistent at infinity for any $\delta > 0$ [details can be found in Alqallaf et al. (2006)].

**6. Concluding remarks.** FDCM assumes that the majority of the cases is clean and follows the underlying model. Robust methods developed for this model exploit the fact that the fraction of clean cases remains constant under affine transformations and concentrate on identifying and downweighting the minority of outlying cases. In fact, the maximal breakdown point of any affine equivariant robust estimator cannot exceed 50%, as shown by Lopuhaä and Rousseeuw (1991). On the other hand, in Section 4 we have shown that the fraction of outlying cells in FICM can drastically change under affine transformations. Consequently, as demonstrated in Section 5, data following FICM and other nonaffine equivariant versions of model (2) can severely upset standard robust procedures, even if the fraction of contaminated cells in the data is quite low.

In practice, both componentwise outliers and structural outliers can occur simultaneously. This situation is considered by the *partially spoiled independent contamination model* (*PSICM*) which assumes that there is a certain probability $\alpha(\epsilon)$ that the case is fully spoiled (as in the FDCM), but otherwise the cells are independently contaminated with probability $\beta(\epsilon)$. Take, for example, $\alpha(\epsilon) = \epsilon/(2-\epsilon)$ and $\beta(\epsilon) = \epsilon/2$.

Similarly, the *partially clean independent contamination model* (*PCICM*) assumes that a case is free of contamination with a certain probability $1 - \alpha(\epsilon)$ (as in the FDCM), but otherwise the different cells are independently contaminated with probability $\beta(\epsilon)$. Two possible choices for the functions $\alpha(\epsilon)$ and $\beta(\epsilon)$ are: (i) $\alpha(\epsilon) = \gamma$ and $\beta(\epsilon) = \epsilon/\gamma$, for some $0 < \gamma < 1$, and (ii) $\alpha(\epsilon) = \beta(\epsilon) = \sqrt{\epsilon}$.

The choices of the functions $\alpha(\epsilon)$ and $\beta(\epsilon)$ in PCICM (i) and (ii) and PSICM are such that the probability of contamination of a single cell is still $P(B_i = 1) = \epsilon$ for all $i$. Therefore, meaningful sensitivity analysis can be performed by letting $\epsilon \to 0$. Another important simplifying feature of these contamination models is that the probability that a case has exactly $k$ contaminated cells does not depend on which are the $k$ contaminated components of the observation. The influence function of multivariate location $M$-estimators under PCICM (i), (ii) and PSICM can be derived similarly as in Section 3 [see Alqallaf et al. (2006) for details]. The influence function under PCICM (i), (ii) turns out to be the same as under FICM and thus is given by (11). The influence function under PSICM becomes

$$\mathrm{IF}(\mu, \mathbf{z}) = \frac{1}{2A_\psi} \left[ \sum_{k=1}^{d} g(H(I_k, \mathbf{z}), \mu_0, \Sigma_0) + g(\Delta_{\mathbf{z}}, \mu_0, \Sigma_0) \right],$$

which is the average of the influence functions under FICM and FDCM. Note that both PCICM and PSICM contain independent componentwise contamination, so the outlier propagation effect occurs in both models. However, the effect will be more devastating in PSICM, where no clean cases are guaranteed. If the fraction of clean cases in PCICM is sufficiently large (at least 50%), then standard affine equivariant robust estimators will show good behavior under this model.

Ideally, robust methods should be resistant against all kind of outliers. However, He and Simpson (1993) showed that the maximal contamination bias of locally linear estimators has to increase with dimension. Moreover, Theorem 1 shows that under FICM the breakdown point of affine equivariant estimators decreases with dimension. These results imply that it is intrinsically difficult to find estimators in high dimensions that are sufficiently robust against all types of outliers. Hence, one has to make a trade-off between several desirable (robustness) properties that cannot all be achieved simultaneously.

Protection against outliers propagation can be achieved by using coordinatewise procedures, such as the (coordinatewise) median, that only operate on one column at the time. Croux et al. (2003) and Maronna and Yohai (2008) use such coordinatewise procedures to construct robust methods for factor models and principal components, respectively. See Liu et al. (2003) for an application to microarray data. However, a well-known weakness of coordinatewise methods is their lack of robustness against structural outliers. This type of outliers can only be handled by robust affine equivariant methods. One possible way to address this trade-off is to apply robust affine equivariant methods to subsets of columns at the time and combine the results. With larger subset sizes, more protection against structural outliers is assured, but less protection against outliers propagation is obtained and vice versa.

## APPENDIX

**A.1. Derivation of (9).** Since $g(H_0, \mu_0, \Sigma) = \mathbf{0}$ for all positive definite matrices $\Sigma$ and elliptically symmetric distributions $H_0$, we have that

$$\frac{\partial}{\partial \Sigma} g(H_0, \mu_0, \Sigma)\Big|_{\Sigma = \Sigma_0} = 0.$$

Hence,

$$\begin{aligned}
(18) \quad & \frac{\partial}{\partial \epsilon} g(H_0, \mu(H(\epsilon, \mathbf{z})), \Sigma(H(\epsilon, \mathbf{z})))\Big|_{\epsilon=0} \\
&= \frac{\partial}{\partial \mu} g(H_0, \mu, \Sigma_0)\Big|_{\mu=\mu_0} \frac{\partial}{\delta \epsilon} \mu(H(\epsilon, \mathbf{z}))\Big|_{\epsilon=0}.
\end{aligned}$$

We also have

$$\begin{aligned}
(19) \quad & \frac{\partial}{\partial \mu} g(H_0, \mu, \Sigma_0)\Big|_{\mu=\mu_0} \\
&= -2 E_{H_0}(\psi'(\mathrm{d}^2(\mathbf{y}, \mu_0, \Sigma_0))(\mathbf{Y} - \mu_0)(\mathbf{Y} - \mu_0)')\Sigma_0^{-1} \\
&\quad - E_{H_0}(\psi(\mathrm{d}^2(\mathbf{Y}, \mu_0, \Sigma_0))\mathbf{I}.
\end{aligned}$$

Let $\mathbf{w} = \Sigma_0^{-1/2}(\mathbf{Y} - \mu_0)$. Then $\mathbf{w}$ has density given by (1) with $\mu_0 = \mathbf{0}$ and $\Sigma_0 = \mathbf{I}$. Since $\mathbf{w}$ has a spherical distribution, it holds that

$$(20) \quad E(\psi'(\|\mathbf{w}\|^2)\mathbf{w}\mathbf{w}') = \frac{1}{d} E(\psi'(\|\mathbf{w}\|^2)\|\mathbf{w}\|^2)\mathbf{I}.$$

From (19) and (20) we get

$$\frac{\partial}{\partial \mu} g(H_0, \mu, \Sigma_0)\Big|_{\mu=\mu_0} = -A_\psi \mathbf{I},$$

where the constant $A_\psi = (2/d)E(\psi'(\|\mathbf{w}\|^2)\|\mathbf{w}\|^2) + E(\psi(\|\mathbf{w}\|^2))$ is independent of $\mu_0$ and $\Sigma_0$. Finally, from (18) we get (9).

**A.2. Proof of Lemma 1.** Suppose that $\mathbf{T}(H)$ is not $\delta$-consistent at infinity for some $\delta > \delta_0$. Then there exists a distribution $G_0$ on $R$, with finite first moment and a sequence of distributions $G_h$ with marginals $G_{hi}(x_i) \leq G_0(x_i - h)$, such that if we call $H_h = (1/2 - \delta)H_0 + (1/2 + \delta)G_h$, then $T_i(H_h) \leq c$ for some $c$, for all $h > 0$. Let

$$H_h^*(\mathbf{x}) = H_h(\sqrt{h}\mathbf{x}) = (1/2 - \delta)H_0(\sqrt{h}\mathbf{x}) + (1/2 + \delta)G_h(\sqrt{h}\mathbf{x}).$$

Then, by scale equivariance of $\mathbf{T}(H)$

(21) $$T_i(H_h^*) \leq c/\sqrt{h} \to 0 \quad \text{as } h \to \infty.$$

Observe that $M_h(\mathbf{x}) = H_0(\sqrt{h}\mathbf{x})$ converges weakly to the point-mass distribution at zero, as $h \to \infty$. Moreover

$$T_i(H_h^*) = (1/2 - \delta) \int x_i w(H_h^*, \mathbf{x}) \, dM_h(\mathbf{x})$$
$$+ (1/2 + \delta) \int x_i w(H_h^*, \mathbf{x}) \, dG_h^*(\mathbf{x}),$$

where $G_h^*(\mathbf{x}) = G_h(\sqrt{h}\mathbf{x})$. Since $w(H_h^*, \mathbf{x}) \leq K$ and $E_{H_0}(|X_i|) < \infty$, we have

$$\left| \int x_i w(H_h^*, \mathbf{x}) \, dM_h(\mathbf{x}) \right| \leq K \int |x_i| \, dM_h(\mathbf{x})$$
$$= K \int \left| \frac{x_i}{\sqrt{h}} \right| dF_0(\mathbf{x}) \to 0 \quad \text{as } h \to \infty.$$

On the other hand, if $A_h = \{\mathbf{x} = (x_1, \ldots, x_d) : x_i \geq h, 1 \leq i \leq d\}$ then $\lim_{h \to \infty} P_{G_h^*}(A_{\sqrt{h}}) = 1$ and therefore $P_{H_h^*}(A_{\sqrt{h}}) \geq 1/2 + \delta$.

Note that by assumptions (iii), $P_{H_h^*}(w(H_h^*, \mathbf{x}) > \eta) > 1/2 - \delta_0$. Set

$$B_h = A_{\sqrt{h}} \cap \{\mathbf{x} : w(H_h^*, \mathbf{x}) > \eta\},$$

then

$$\lim_{h \to \infty} P_{H_h^*}(B_h) \geq \lim_{h \to \infty} P_{H_h^*}(\{x : w(H_h^*, \mathbf{x}) > \eta\}) - \lim_{h \to \infty} P_{H_h^*}(A_{\sqrt{h}}^c)$$
$$\geq 1/2 - \delta_0 - (1/2 - \delta) = \delta - \delta_0 > 0.$$

Since $\lim_{h \to \infty} P_{M_h}(B_h) = 0$, we have

$$0 < \delta - \delta_0 \leq \lim_{h \to \infty} P_{F_h^*}(B_h)$$
$$= (1/2 - \delta) \lim_{h \to \infty} P_{M_h}(B_h) + (1/2 + \delta) \lim_{h \to \infty} P_{G_h^*}(B_h)$$
$$= (1/2 + \delta) \lim_{h \to \infty} P_{G_h^*}(B_h).$$

Therefore, $\lim_{h\to\infty} P_{G_h^*}(B_h) \geq \gamma = (\delta - \delta_0)/(1/2 + \delta)$. Then

$$
\begin{aligned}
\lim_{h\to\infty} &\int x_i w(H_h^*, \mathbf{x}) \, dG_h^*(\mathbf{x}) \\
&\geq \sqrt{h}\eta \lim_{h\to\infty} \int_{B_h} dG_h^*(\mathbf{x}) + \lim_{h\to\infty} \int_{x_i<0} x_i w(H_h^*, \mathbf{x}) \, dG_h^*(\mathbf{x}) \\
&= \sqrt{h}\eta \lim_{h\to\infty} \int_{B_h} dG_h^*(\mathbf{x}) + \lim_{h\to\infty} \int_{x_i<0} \frac{x_i}{\sqrt{h}} w\left(H_h^*, \frac{\mathbf{x}}{\sqrt{h}}\right) dG_h(\mathbf{x}) \\
&\geq \sqrt{h}\eta \lim_{h\to\infty} \int_{B_h} dG_h^*(\mathbf{x}) + K \lim_{h\to\infty} \int_{x_i<0} \frac{x_i}{\sqrt{h}} \, dG_h(\mathbf{x}).
\end{aligned}
\tag{22}
$$

Now, regarding the first term we have of the right-hand side of (22) we get

$$
\eta \lim_{h\to\infty} \sqrt{h} \int_{B_h} dG_h^*(\mathbf{x}) = \eta \lim_{h\to\infty} \sqrt{h} P_{G_h^*}(B_h) \geq \gamma\eta \lim_{h\to\infty} \sqrt{h} = \infty.
\tag{23}
$$

Regarding the second term, first note that the distribution of $x_i (x_i < 0)$ under $G_h(\mathbf{x})$ is also stochastically larger than the corresponding distribution under $G_0(x - h)$. Then, by the change of variable $y = x - h$, and using this stochastic inequality we have

$$
\begin{aligned}
\int_{x_i<0} x_i \, dG_h(\mathbf{x}) &\geq \int_{x_i<0} x_i \, dG_0(x_i - h) = \int_{y_i<-h} (y_i + h) \, dG_0(y_i) \\
&= \int_{y_i<-h} y_i \, dG_0(y_i) + h \int_{y_i<-h} dG_0(y_i).
\end{aligned}
\tag{24}
$$

The first term is uniformly bounded as follows

$$
\left|\int_{y_i<-h} y_i \, dG_0(y_i)\right| \leq \int_{y_i<-h} |y_i| \, dG_0(y_i) \leq \int_{-\infty}^{\infty} |y_i| \, dG_0(y_i) < \infty.
\tag{25}
$$

The second term tends to zero because $G_0$ has finite first moments and so

$$
h \int_{y_i<-h} dG_0(y_i) = h P_{G_0}(y_i < -h) \to 0.
\tag{26}
$$

By (23) the first term in (22) tends to $+\infty$. By (24), (25) and (26) the second term in (22) is uniformly bounded. Therefore,

$$
\lim_{h\to\infty} T_i(H_h^*) = \lim_{h\to\infty} \int x_i w(H_h^*, \mathbf{x}) \, dG_h^*(\mathbf{x}) = +\infty,
$$

contradicting (21).

**A.3. Proof of Lemma 2.** The weights $w(\mathbf{x}, H)$ in (16) are clearly nonnegative. Moreover, $u(t)$ is bounded because by assumption (A2), $u(t) \leq u(0) =$ $\psi'(0) = \kappa < \infty$. For any $0 < t < 1$, let $B_t = \{\mathbf{x} : \rho(\mathrm{d}(\mathbf{x}, \mathbf{T}(H), S(H))/s_0) \leq t\}$. It follows from (13) and (14) that

$$\tfrac{1}{2} \geq \int_{(B_t)^c} \rho(\mathrm{d}(\mathbf{x}, \mathbf{T}(H), S(H))/s_0)) \, dH(\mathbf{x}) \geq (1 - P(B_t))t$$

and so

$$P(B_t) \geq 1 - \frac{1}{2t}. \tag{27}$$

Let $t_0$ be such that

$$1 - \frac{1}{2t_0} = \frac{1}{2} - \delta_0, \tag{28}$$

that is, $t_0 = 1/(1 + 2\delta_0)$. Note that $0 < \delta_0 \leq 1/2$ implies that $1/2 \leq t_0 < 1$. Moreover, combining (27) and (28) yields

$$P(B_{t_0}) \geq \tfrac{1}{2} - \delta_0. \tag{29}$$

Let $r_0 = \rho^{-1}(t_0)$, then we can write $B_{t_0} = \{\mathbf{x} : \mathrm{d}(\mathbf{x}, \mathbf{T}(H), S(H))/s_0 \leq r_0\}$. By monotonicity of $u(t)$, for all $\mathbf{x}$ in $B_{t_0}, u(d(\mathbf{x}, \mathbf{T}(H), S(H))/s_0) \geq u(r_0)$. Put $\zeta = u(r_0)$, since $1/2 \leq t_0 < 1$, using the assumption that $\psi(u) > 0$ when $\rho(u) < 1$ we have $\zeta > 0$. Then we can write $u(\mathrm{d}(\mathbf{x}, \mathbf{T}(H), S(H))/s_0) \geq \zeta$, for $\mathbf{x} \in B_{t_0}$ and then $\kappa \geq E(u(\mathrm{d}(\mathbf{X}, \mathbf{T}(H), S(H))/s_0)) \geq \zeta P(B_{t_0}) \geq \zeta/4$. Therefore, $w(\mathbf{x}, H) \leq 4\kappa/\zeta$ for all $\mathbf{x}$ and $w(\mathbf{x}, H) > \zeta/\kappa$ for $\mathbf{x} \in B_{t_0}$. Together with (29) this means that assumptions (i)–(iii) of Lemma 1 are satisfied.

F. Alqallaf
Department of Statistics
and Operations Research
Faculty of Science
Kuwait University
P.O. Box 5969 Safat-13060
Kuwait
E-mail: fatemah@kuc01.kuniv.edu.kw

S. Van Aelst
Department of Applied Mathematics
and Computer Science
Ghent University
Krijgslaan 281 S9
B-9000 Gent
Belgium
E-mail: Stefan.VanAelst@UGent.be

V. J. Yohai
Department of Mathematics
University of Buenos Aires
Ciudad Universitaria, Pabellón 1
1426 Buenos Aires
Argentina
E-mail: vyohai@dm.uba.ar

R. H. Zamar
Department of Statistics
University of British Columbia
6356 Agricultural Road
Vancouver, British Columbia
Canada V6T 1Z2
E-mail: ruben@stat.ubc.ca